\documentclass[12pt]{amsart}
\usepackage{latexsym,amsfonts,amsthm,amsmath,amscd,amssymb}

\newtheorem{theorem}{Theorem}

\newtheorem{lemma}[theorem]{Lemma}

\begin{document}

\title{On the Gauss Map with Vanishing Biharmonic stress-energy tensor}

\author{Wei Zhang}

\thanks{Mathematics Classification Primary(2000): 58E20.
\indent Thank Prof. Dong and Prof. Ji for helpful discussion.}

\maketitle
\begin{abstract}
We study the biharmonic stress-energy tensor $S_2$ of Gauss map.
Adding few assumptions, the Gauss map with vanishing $S_2$ would be
harmonic.
\end{abstract}

\section{Introduction}
Let $\phi: (M,g)\rightarrow(N,h)$ be a smooth map between two
Riemannian manifolds. Assume M compact and define the energy of
$\phi$ to be:

$$E(\phi)=\frac{1}{2}\int_M |d\phi|^2 v_g$$

Call a map harmonic if it is a critical point of E, and this is
characterized by:

$$\tau(\phi)=trace \nabla d\phi=0$$
where $\tau$ is called the tension field. The vanishing of this
field is used to define the harmonic map in noncompact case.

As a natural generalization, the biharmonic map is the critical
point of the bienergy(see details in \cite{MO}):

$$E_2(\phi)=\frac{1}{2}\int_M |\tau(\phi)|^2 v_g$$
it's associated bitension field is:

$$\tau_2(\phi)=-\Delta^\phi \tau(\phi)-trace
R^N(d\phi,\tau(\phi))d\phi$$

Described by Hilbert \cite{Hi}, the stress-energy tensor associates
to a variational problem is a symmetric 2-covariant tensor S
conservative at critical points, i.e. $div S=0$.

In the context of harmonic maps, the stress-energy tensor was
studied  by Baird and Eells in details(\cite{BE}):

$$S=\frac{1}{2}|d\phi|^2 g-\phi^* h$$
s.t $div S=-<\tau(\phi),d\phi>$.

The biharmonic stress-energy tensor $S_2$ corresponding to
$E_2(\phi)$ was introduced by Jiang in \cite{J} and studied by
Loubeau et.al \cite{LMO1}:

\begin{equation}
\begin{split}
 S_2(X,Y)=&\frac{1}{2}|\tau(\phi)|^2<X,Y>+<d\phi,\nabla
\tau(\phi)><X,Y>\\
&-<d\phi(X),\nabla_Y \tau(\phi)>-<d\phi(Y),\nabla_X \tau(\phi)>
\end{split}
\end{equation}

It can be easily drawn out: $S=0\Rightarrow \tau=0 \Rightarrow S_2=0
\Rightarrow \tau_2=0$.

A natural question is when the reverse be true. For example, when a
biharmonic map be harmonic? Readers could confer \cite{MO} for the
results in this direction.

In this article, we focus on how would $S_2=0$ imply $\tau=0$. Jiang
had proved:

\begin{theorem}[\cite{J}]
A map $\phi$: $(M,g)\rightarrow (N,h)$, $m \neq 4$, with $S_2=0$, M
compact and orientable, is harmonic
\end{theorem}

\begin{proof}
Trace of $S_2$ is

$$0=trace S_2=\frac{m}{2}|\tau(\phi)|^2+m<d\phi,\nabla
\tau(\phi)>-2<d\phi,\nabla \tau(\phi)>$$ %
and integrating over M:

$$0=\frac{4-m}{2}\int_M |\tau(\phi)|^2 v_g$$
hence, $\phi$ is harmonic when $m \neq 4$.
\end{proof}

And

\begin{theorem}[\cite{J}]
A non-minimal Riemannian immersion $\phi$: $(M^4,g)\rightarrow
(N,h)$ saisfies $S_2=0$ if and only if it is pseudo-umbilical.
\end{theorem}

The 4-folds have somehow exotic behaviors. At another hand, if the
map $\phi$ arise from a submanifold's Gauss map, under few
assumptions, the map would regain part of "rigidity", i.e. $S_2=0$
implying harmonic. It is the theme of this paper:

\begin{theorem}
$M^4$ is a compact pseudo-umbilical submanifold of $R^n$, if its
Gauss map has vanishing $S_2$, then $M^4$ has constant mean
curvature.
\end{theorem}

And

\begin{theorem}
$M^4$ is a compact analytic hypersurface of $R^5$, if it is strictly
convex and its Gauss map has vanishing $S_2$, then M is a
hypersphere, i.e the Gauss map is identity.
\end{theorem}

\section{preliminary}

Form now on, $M^m$ always denotes an oriented submanifold of $R^n$,
and $G$ shorts for the Gauss map.

In order to study the biharmonic stress-energy tensor of Gauss map,
we have to understand $G^* T (G(n,m))$ well, especially the
connection on it.

In \cite{RV}, Ruh and Vilms had shown: the pull back of the tangent
bundle of Grassmannian via the Gauss map is isomorphic to
$T^*(M)\otimes N(M)$, i.e. $T(M)\otimes N(M)$ after the musician
transformation, where $T(M)$, $T^*(M)$, $N(M)$ are the tangent,
cotangent and normal bundle respectively.

There is a more explicit way to see this(\cite{LMO2}):

Choose $\{e_i\}_{i=1}^m$ an oriented geodesic basis centered around
$p \in M$. In the neighborhood $U \ni p$, the Gauss map can be
written as:

$$G(q)=e_1(q)\wedge \cdots \wedge e_m(q), \forall q \in U$$

Since

$$dG_q(e_i)=\sum_{j=1}^m e_1(q)\wedge \cdots \wedge e_j(q)
\wedge(\nabla^{R^n}_{e_i}e_j)(q)\wedge e_{j+1}(q)\wedge \cdots
\wedge e_m(q),$$ %
restricting at p, we have:

$$dG_p(e_i)=\sum_{j=1}^m e_1(p)\wedge \cdots \wedge e_j(p)
\wedge B_p(e_i,e_j)\wedge e_{j+1}(p)\wedge \cdots \wedge e_m(p)$$
where B is the second fundamental form of M, taking value in N(M).
Now $dG_p(e_i)$ can be identified with $\sum_j e_j^*(p) \otimes
B_p(e_i,e_j)$, i.e. $\sum_j e_j(p) \otimes B_p(e_i,e_j)$. Thus the
bundles are isomorphic.

For latter utility, we'd better explain the relationship between
above invariant method and the moving frames method. See \cite{JX}
or \cite{X}, complete $\{e_i(p)\}_{i=1}^m$ into an orthonormal basis
$\{e_\alpha(p)\}_{\alpha=1}^n$ of $R^n$. $\{w_\alpha\}$ is the dual
frame. The Riemannian connection on $R^n$ is uniquely determined by
the equation:

$$d w_\alpha=w_{\alpha \beta} \wedge w_\beta$$
$$w_{\alpha \beta}+w_{\beta \alpha}=0$$

The canonical metric on G(m,n) is:

\begin{equation}\label{metric}
d s^2=\sum_{i,a} w_{ia}^2
\end{equation}
where $a=m+1, \cdots, n$.

$\{w_{ia}\}$ can be thought as the dual frame of $\{e_i \otimes
e_a\}$, which means nothing but $e_j \otimes e_a$ is orthonormal
basis of $T_{G(p)}G(n,m)$ with respect to the canonical metric on
Grassmannian, i.e.

$$g_{can}(dG_p(e_i),dG_p(e_k))=\sum_j<B_p(e_i,e_j),B_p(e_k,e_j)>$$

The connection due to metric given by equation ~(\ref{metric}) is:

$$w_{ia jb}=\delta_{ab}w_{ij}+\delta_{ij}w_{ab}$$

Pull back these forms to M, we have:

\begin{lemma}[\cite{X}]
The connection on $G^* T (G(n,m))$ is the connection $\nabla^M
\otimes \nabla^\bot$ on $T(M) \otimes N(M)$
\end{lemma}

\section{Proof of main theorem}

In \cite{RV}, it had been showed that the tension field of the Gauss
map is identical with $\nabla H$, $\nabla_{e_i}H \otimes e_i$ in our
setting, where H is the mean curvature.

Let G be the Gauss map of $M^4$. Take the trace of its biharmonic
stress-energy tensor $S_2$:

\begin{equation}\label{trace}
\begin{split}
\frac{1}{2}trace S_2=&|\tau(G)|^2+<\nabla \tau(G),dG> %
=\nabla_{e_i}<\tau (G), dG(e_i)>\\
=&\nabla_{e_i}<\nabla_{e_j}H \otimes e_j,B_{ik} \otimes e_k> %
=\nabla_{e_i}<\nabla_{e_j}H, B_{ij}>\\
=&\nabla_{e_i} \nabla_{e_j} <H, B_{ij}>-\nabla_{e_i}<H,
\nabla_{e_i}B_{jj}>
\end{split}
\end{equation}

If $M^4$ is pseudo-umbilical and  $S_2$ vanishes, equation
~(\ref{trace}) become:

$$0=\frac{1}{4}\nabla_{e_i} \nabla_{e_i}|H|^2-\frac{1}{2}\nabla_{e_i}
\nabla_{e_i}|H|^2 %
=-\frac{1}{4}\Delta|H|^2$$

$|H|^2$ is a harmonic function, by the maximum principle, it must be
constant. This ends the proof of the theorem 3.

\noindent {\bf Remark:} By theorem 2, we know that in the case m=4,
if an embedding and its Gauss map are both with vanishing biharmonic
stress-energy tensor, then the submanifold has constant mean
curvature.

 For the sequel, we need a reformulation of $S_2=0$:

\begin{lemma}[\cite{LMO1}]
Let $\phi: (M^4,g)\rightarrow(N,h)$, then $S_2=0$ if and only if

\begin{equation}\label{reformS}
\begin{split}
\frac{1}{2}&|\tau(\phi)|^2<X,Y>+<d\phi(X),\nabla_Y \tau(\phi)>\\
&+<d\phi(Y),\nabla_X \tau(\phi)>=0
\end{split}
\end{equation}
$\forall X,Y$
\end{lemma}

Now, we are in the position to prove theorem 4.

\begin{proof}

Denote h the scalar value of H, by equation ~(\ref{reformS}) we
have:

\begin{equation}\label{invEq}
\begin{split}
0=&\frac{1}{2}|gradh|^2<X,Y>+<B(X,e_j)e_j,\nabla_Y gradh>\\ %
&+<B(Y,e_j)e_j,\nabla_X gradh>\\ %
=&\frac{1}{2}|gradh|^2<X,Y>+B(X,\nabla_Y gradh)+B(Y,\nabla_X gradh)
\end{split}
\end{equation}

For M is compact, assume h achieves its maximum at point p. Choose
local patch around p $\{U, x_i\}$, s.t $\{ \frac{\partial}{\partial
x_i}\}$ is orthnormal at p and diagonalizing the second fundamental
form. Replace X, Y with $\frac{\partial}{\partial x_i}$,
$\frac{\partial}{\partial x_j}$, write equation ~(\ref{invEq}) in
local coordinates:

\begin{equation}\label{localEq}
\frac{\partial h}{\partial x_l \partial x_j}B_{il} %
+\frac{\partial h}{\partial x_l \partial x_i}B_{jl} %
=-\frac{\partial h}{\partial x_l}\Gamma_{jl}^k B_{ik} %
-\frac{\partial h}{\partial x_l}\Gamma_{il}^k B_{jk} %
-\frac{1}{2} \frac{\partial h}{\partial x_l} g_{lk} \frac{\partial
h}{\partial x_k} g_{ij}
\end{equation}

If B has eigenvalues $\{\lambda_i\}$, equation ~(\ref{localEq}) at
point p is:

$$\frac{\partial h}{\partial x_i \partial x_j}(\lambda_i+\lambda_j) %
=-\frac{1}{2} \frac{\partial h}{\partial x_l} g_{lk} \frac{\partial
h}{\partial x_k} g_{ij} %
-\frac{\partial h}{\partial x_l}\Gamma_{jl}^k B_{ik} %
-\frac{\partial h}{\partial x_l}\Gamma_{il}^k B_{jk} %
$$

M is analytic, so we treat every thing in the analytic category.

Expend h as a convergent polynomial series in neighborhood of p:
$$h(x_1,\cdots,x_4)=c+h_{ij}x_ix_j+h_{ijk}x_ix_jx_k+\cdots$$
where $h_{ij}$ means $h_{ij}(p)$ in effect.

It has no one order terms.

Expend B either, we have:
\begin{equation}\label{ComOdr}
\begin{split}
&\frac{\partial h}{\partial x_i \partial x_j}(\lambda_i+\lambda_j) %
+ O(x)F(h_{kl}) \\ %
 =&-\frac{1}{2} \frac{\partial h}{\partial x_l} g_{lk} \frac{\partial
h}{\partial x_k} g_{ij} %
-\frac{\partial h}{\partial x_l}\Gamma_{jl}^k B_{ik} %
-\frac{\partial h}{\partial x_l}\Gamma_{il}^k B_{jk} %
\end{split}
\end{equation}
where F is a linear combination if $h_{ij}$.

Since:

$$\frac{\partial h}{\partial x_i \partial x_j} %
=h_{ij}+\sum_k h_{ijk}x_k+O(x^2)$$

$$\frac{\partial h}{\partial x_i}  %
=\sum_k h_{ik}x_k+O(x^2)$$

Comparing the lowest order terms in two sides of equation
~(\ref{ComOdr}).

For $\lambda_i>0$, the left hand side may contains zero terms, while
the left hand side has order no less than one, so $h_{ij}$ must
vanish, $\forall i,j$.

Using the boot-strap argument, it is easy to see all the derivative
of h at p must be zero.

Thus h must be constant for the anality. This means the Gauss map is
harmonic.

In fact we have more stronger conclusion. Hsiung had shown in
\cite{Hs} that a strictly convex hypersurfaces with constant mean
curvature must be hypersphere.

\end{proof}

\newpage

\noindent Wei Zhang

\

\noindent School of Mathematical Sciences

\noindent Fudan University

\noindent Shanghai, 200433, P. R.China

\

\noindent Email address: 032018009@fudan.edu.cn
\end{document}